\documentclass[10pt,a4paper]{article}
\usepackage{amssymb}
\usepackage{latexsym,bm}
\usepackage{graphicx}
\usepackage{mathrsfs,amscd,amssymb,amsthm,amsmath,bm,graphicx,psfrag,subfigure,url,xcolor,ifpdf}
\usepackage{amsmath}
\usepackage{mathrsfs}

\date{}
\begin{document}
\title{ Extremal polygonal chains with respect to the Kirchhoff index  \footnote{E-mail addresses:
{\tt mathqima@163.com}(Q. Ma).}}
\author{\hskip -10mm  Qi Ma\thanks{Corresponding author.}\\
{\hskip -10mm \small Ministry of Basic Education, Sanda University, Shanghai 200241, China}}\maketitle
\begin{abstract}
The Kirchhoff index is defined as the sum of resistance distances between all pairs of vertices in a graph. This index is a critical parameter for measuring graph structures.  In this paper, we characterize polygonal chains with the minimum Kirchhoff index, and characterize even (odd) polygonal chains with the maximum Kirchhoff index, which extends the results of \cite{45}, \cite{14} and \cite{2,13,44} to a more general case.
\end{abstract}


\vskip 8mm

\section{ Introduction}
We consider finite simple graphs. The {\it order} of a graph is the number of vertices, and the {\it size} is the number of edges. We denote by $V(G)$ and $E(G)$ the vertex set and edge set of a graph $G$ respectively. Denote by $d_G(u,v)$ the {\it distance} between two vertices $u$ and $v$ in $G.$ We use $G+H$ to denote the {\it disjoint union} of graphs $G$ and $H.$ The {\it join} of graphs $G$ and $H,$ written $G\vee H,$ is the graph obtained from $G+H$ by adding all the edges between $G$ and $H.$ For the terminology and notations not defined here, we follow the book \cite{7}.

Molecules can be modeled by graphs with vertices for atoms and edges for atomic bonds. The topological indices of a graph can provide some information on the chemical properties of the corresponding molecule. One of the most famous indices is the {\it Wiener index} \cite{8} which is defined as
$$
 W(G)=\sum_{\{u,v\}\subseteq V(G)}d_G(u,v).
$$
Wiener \cite{8} used it to study the boiling point of paraffin. Many chemical properties of molecules are related to the Wiener index \cite{15,21,22,16}. Based on the electronic network theory, Klein and Randi\'{c} \cite{9} proposed the concept of {\it resistance distance} in 1993. The resistance distance $r_G(u,v)$ between vertices $u$ and $v$ of a connected graph $G$ is computed as the effective resistance between vertices $u$ and $v$ in the corresponding electrical network constructed from $G$ by replacing each edge of $G$ with an unit resistor. This novel parameter is in fact intrinsic to the graph and has some nice interpretations and applications in chemistry; see \cite{20,18,17,19} for more detailed information. Similar to the Wiener index, Klein and Randi\'{c} \cite{9} defined the {\it Kirchhoff index} $K\!f(G)$ of $G$ as the sum of the resistance distances between all pairs of vertices in $G,$ i.e.
$$
 K\!f(G)=\sum_{\{u,v\}\subseteq V(G)}r_G(u,v).
$$

Recently, more and more attention has been paid to the Kirchhoff index. An important research direction about the Kirchhoff index is determining the graphs with the maximal or minimal Kirchhoff index in a given class of graphs. Up to now, the Kirchhoff index has been given for many graphs, such as silicate networks \cite{32}, ladder-like graphs \cite{24}, ladder graphs \cite{23}, linear pentagonal chains \cite{30}, linear hexagonal chains \cite{25}, linear crossed hexagonal chains \cite{26}, linear crossed octagonal graphs \cite{10}, phenylenes chians \cite{27}, cyclic phenylenes \cite{29}, M\"{o}bius phenylenes chains and cylinder phenylenes chains \cite{28}, linear octagonal chains \cite{12}, M\"{o}bius/cylinder octagonal chains \cite{11}, random cyclooctatetraene and spiro chains \cite{36}, random polyphenyl and spiro chains \cite{31}. Some other topics on the Kirchhoff index of a graph may be referred to \cite{33,34,35,46,38} and references therein.

A planar graph is said to be an {\it outerplane graph} if all the vertices lie on the boundary of the exterior face. Let $P_n$ be a $2$-connected outerplane graph and satisfying the following three conditions. (1) $P_n$ has $n$ interior faces and the length of each interior face is at least four. (2) Any two interior faces share exactly one edge , or disjoint. (3) Any vertex has degree $2$ or $3$.  Note that each interior face of $P_n$ is a polygon, so we called $P_n$ a {\it polygonal chain}. If all the polygons of $P_n$ have an even (odd) size, we call $P_n$ an {\it even (odd) polygonal chain}. If all the polygons of $P_n$ have the same size $l,$ we call $P_n$ an {\it $\ell$-polygonal chain}. Recently, a lot of attention has been paid to the  polygonal chain. Zhang and Jiang \cite{40} studied the continuous forcing spectra of even polygonal chains. Chen and Li \cite{15} determined the expected values of Wiener indices in random even polycyclic chains. The extremal polygonal chains on $k$-matchings are characterized by Cao and Zhang \cite{43}. Wei and Shiu \cite{39} studied the Wiener indices of a random $\ell$-polygonal chain and its asymptotic behavior, which covering some previous results for special random chains. Zhu, Fang and Geng \cite{41} determined the Gutman and Schultz indices in the random $\ell$-polygonal chains. Zhu and Geng \cite{42} determined the multiplicative degree-Kirchhoff index in the random polygonal chains. Li and Wang \cite{37} established an explicit analytical expressions for the expected values of the (degree)-Kirchhoff index in random $\ell$-polygonal chains.

Inspired by these excellent works, we naturally consider the extremal polygonal chain $P_n$ with respect to Kirchhoff index. We first introduce some notations. We denote the polygons of $P_n$ by $H_1,H_2,\ldots,H_n$ respectively such that $H_i$ is adjacent to $H_{i+1},(1\leq i\leq n-1),$ and let $k_i$ be the length of $H_i.$ Denote by $Q_{n}$ a ladder graph with $n$ squares and denote by $S_i   (i=1,\ldots,n)$ the $i$-th square of $Q_{n}$. Note that a polygonal chain $P_n$ can be obtained from $Q_{n}$ by adding $k_i-4$ vertices to the $i$-th square $S_i$ of $Q_{n}$. We have $k_i-3$ ways to add these $k_i-4$ new vertices to $S_{i}.$ That is, we can add $0$ (resp. $1,\ldots, k_i-4$) vertices to the top edge of $S_{i}$ and the remaining vertices to the bottom edge of $S_{i}.$ Obviously, adding a vertex to the top edge of $S_1(S_n)$ or bottom edge of $S_1(S_n)$ does not affect the structure of the graph $P_n$. For convenience, we always suppose that $H_1$ and $H_n$ are obtained by adding $k_i-4$ to the bottom edge of $S_{1}$ and $S_n.$ For the  polygon $H_i (2\leq i\leq n-1)$, we give a number $w_i=0\,$(resp. $1,2,\ldots, k_i-4$) to $H_i$ if $H_i$ is obtained by adding $w_i$ vertices to the top edge of $S_i$. In this viewpoint, we are able to represent $P_n$ by a vector $w=(w_2,w_3,\ldots ,w_{n-1})$ with $w_i\in \{0,1,2,\ldots,k_i-4\}$. In the following, we always denote a polygonal chain with $n$ polygons by $P(w)$ such that $w$ is an $(n-2)$-tuple of $\{0,1,2,\ldots,k_i-4\}$.

Next we introduce some special polygon chains. Let $2\leq i\leq n-1,$ a {\it kink} in a polygonal chain is a polygon with $w_i=0$ or $k_i-4,$ and a polygonal chain $P(w)$ with $w_i=0$ or $k_i-4$ for all $i$ is called a {\it ``all-kink" chain}. The polygonal chain $P(\underbrace{0,0,\ldots,0}_{n-2})$ (isomorphic to $P(\underbrace{k_2-4,k_3-4,\ldots,k_{n-1}-4}_{n-2})$) is called a {\it helicene polygonal chain}, denoted by $HP_n$. The $7$-polygonal chain $HP_5$ is shown in Fig. 1. If a polygon chain has $w_i=\lfloor\frac{k_i-4}{2}\rfloor$ or $\lceil\frac{k_i-4}{2}\rceil,2\leq i\leq n-2,$ we call this polygon chain a {\it linear polygon chain} and denote it by $LP_n.$
\begin{figure}[!ht]
  \centering
 \includegraphics[width=110mm]{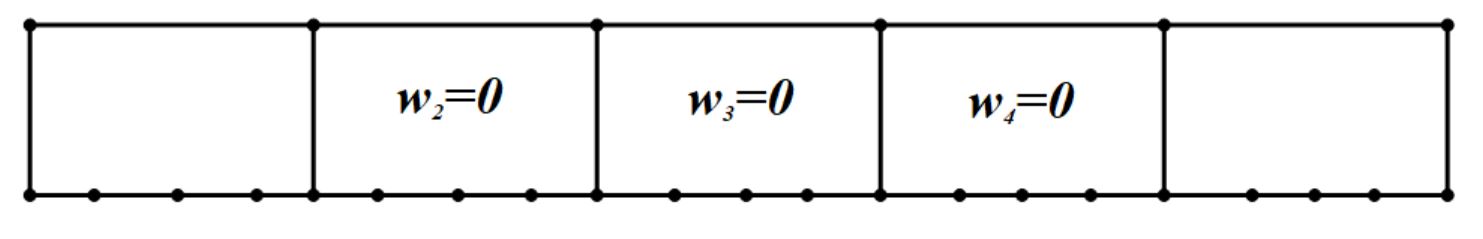}
 \caption{The graph $HP_5$ with all polygons have size $7$.}
\end{figure}

Very recently, Yang and Sun \cite{13} characterized hexagonal chains ($k_i=6,1\le i\le n$) with minimum Kirchhoff indices. They \cite{14} also characterized pentagonal chains ($k_i=5,1\le i\le n$) with minimum Kirchhoff indices. Generalizing these two results, Liu and You \cite{45} characterized $\ell$-polygonal chains with minimum Kirchhoff indices. A phenylene chain is a polygonal chain with $k_i=6$ for odd $i$ and $k_i=4$ for even $i.$ For phenylene chains, Yang and Wang \cite{44} conjecture that the helicene chain has the minimum Kirchhoff index. Motivated by previous researchs, we characterize polygonal chains with minimum Kirchhoff indices. Our main results are as follows:

{\bf Theorem 1.} {\it Among all polygonal chains, the helicene polygonal chain is the unique chain with minimum Kirchhoff index. }

In \cite{2}, Yang and Klein characterized hexagonal chains ($k_i=6,1\le i\le n$) with maximum Kirchhoff index. Sun and Yang \cite{14} characterized pentagonal chains ($k_i=5,1\le i\le n$) with maximum Kirchhoff indices. Yang and Wang \cite{44} characterized phenylene chains  with maximum Kirchhoff indices.  Liu and You \cite{45} characterized $\ell$-polygonal chains with maximum Kirchhoff indices. Generalizing the above results, we have the following theorem and corollaries.

{\bf Theorem 2.} {\it Among all polygonal chains, the maximum Kirchhoff index is obtained only when the polygonal chain is a linear polygonal chain. }

By Theorem 2, we have the following corollaries which characterizes even (odd) polygonal chains with the maximum Kirchhoff indices.

{\bf Corollary 3.} {\it Let $k_i$ denote the length of the $i$-th polygon in a polygonal chain. Then among all even polygonal chains, the polygonal chain $P(\frac{k_2-4}{2},\frac{k_3-4}{2},\ldots,\frac{k_{n-1}-4}{2})$ has the maximum Kirchhoff index .}

{\bf Corollary 4.} {\it Let $k_i$ denote the length of the $i$-th polygon in a polygonal chain. Then among all odd polygonal chains, the polygonal chain $P(\frac{k_2-5}{2},\frac{k_3-3}{2},\frac{k_{4}-5}{2},\frac{k_5-3}{2},\ldots)$ has the maximum Kirchhoff index. }

\section{Preliminaries}

To prove our main results, we will need the following definitions and lemmas.

{\bf Definition 1.} (Series Transformation) Let $x, y$ and $z$ be nodes in a graph where $y$ is adjacent to only $x$ and $z$. Moreover, let $R_1$ equal the resistance between $x$ and $y$ and $R_2$ equal the resistance between $y$ and $z$. A series transformation transforms this graph by deleting $y$ and setting the resistance between $x$ and $z$ equal to $R_1 + R_2$.

{\bf Definition 2.} (Parallel Transformation) Let $x$ and $y$ be nodes in a multi-edged graph where $e_1$ and $e_2$ are two edges between $x$ and $y$ with resistances $R_1$ and $R_2$, respectively. A parallel transformation transforms the graph by deleting edges $e_1$ and $e_2$ and adding a new edge between $x$ and $y$ with edge resistance $r=(\frac{1}{R_1}+\frac{1}{R_2})^{-1}.$

A $\Delta$-$Y$ transformation is a mathematical technique to convert resistors in a triangle formation to an equivalent system of three resistors in a $``Y"$ format as shown in Fig. 2. We formalize this transformation below.
\begin{figure}[!ht]
\centering
  \includegraphics[width=80mm]{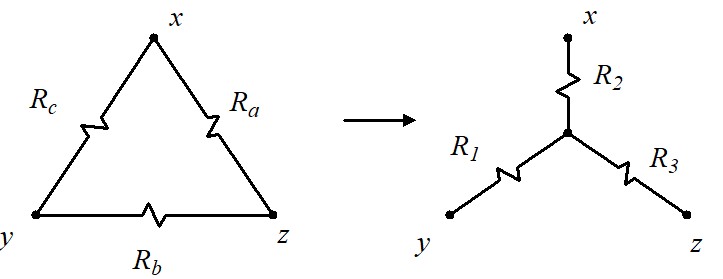}
  \caption{$\Delta$ and $Y$ circuits with vertices labelled as in Definition 3.}
\end{figure}

{\bf Definition 3.} ($\Delta$-$Y$ Transformation) Let $x,y,z$ be nodes and $R_a,R_b$ and $R_c$ be given resistances as shown in Fig. 2. The transformed circuit in the $``Y"$ format as shown in Fig. 2 has the following resistances:
\begin{eqnarray*}
  R_1=\frac{R_bR_c}{R_a+R_b+R_c}, \ \  R_2=\frac{R_aR_c}{R_a+R_b+R_c}, \ \
  R_3=\frac{R_aR_b}{R_a+R_b+R_c}.
\end{eqnarray*}

{\bf Lemma 5.} (Stevenson \cite{6}) {\it Series transformations, parallel transformations, and $\Delta-Y$ transformations yield equivalent circuits.}

{\bf Definition 4.} In the following, we will use $r_{G}(x)$ to denote the sum of resistance distances between $x$ and each other vertex of $G.$ More precisely,
$
r_{P_n}(x)=\sum_{y\in V(P_n)}r_{G}(x,y).
$

Let $A$ and $B$ be two vertex-disjoint graphs, $u,v\in V(A)$ and $u\neq v,$ $x,y\in V(B)$ and $x\neq y$. Let $S$ be the graph obtained from $A$ and $B$ by connecting $u$ with $x$, and $v$ with $y$ as shown in Fig. 3. Let $T$ be the graph obtained from $S$ by deleting edges $\{ux,\,vy\}$ and adding edges $\{uy,\,vx\};$ See Fig. 3. Then we call $S$ and $T$ are $S,T$-isomers. The concept of $S,T$-isomers was introduced by Polansky and Zander \cite{5} in 1982. From then on, a lot of research \cite{1,3,4,2} has been devoted to the study of topological properties of $S,T$-isomers.
\begin{figure}[!ht]
\centering
  \includegraphics[width=40mm]{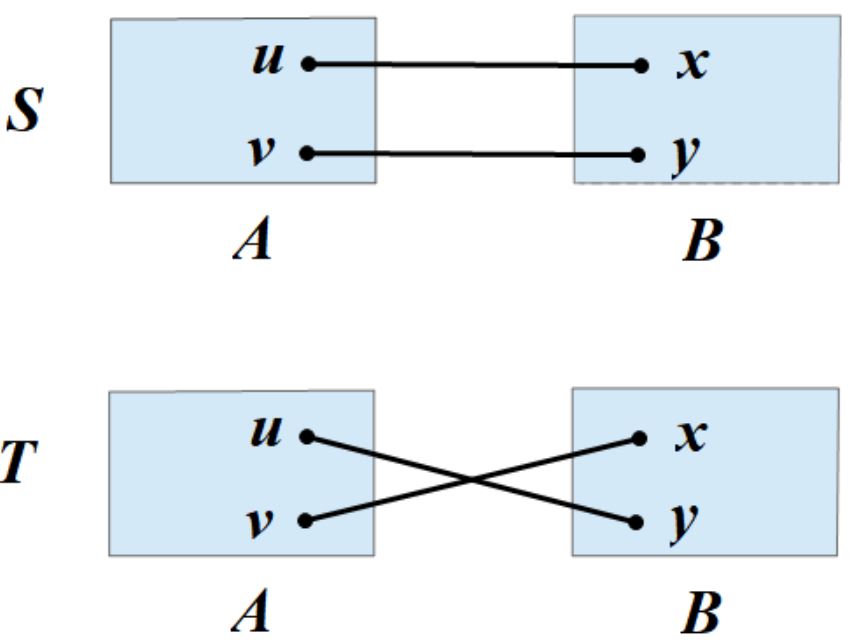}
  \caption{The structure of the graphs $S$ and $T$ and the labeling of their vertices.}
\end{figure}

Yang and Klein \cite{2} obtained the comparison theorem on the Kirchhoff index of $S,T$-isomers, which plays an essential role in the characterization of extremal polygonal chains.

{\bf Lemma 6.} (Yang and Klein \cite{2}) {\it Let $S,T,A,B,u,v,x,y$ be defined as shown in Fig. 3. Then}
$$
K\!f(S)-K\!f(T)=\frac{[r_A(u)-r_A(v)][r_B(y)-r_B(x)]}{r_A(u,v)+r_B(x,y)+2}.
$$

{\bf Lemma 7.} (Klein and Randi\'{c} \cite{3}) {\it The resistance function on a graph is a distance function. Thus for any vertices $x,y,z\in V(G)$, we have\\
(1) $r_G(x,y)\ge 0$,\\
(2) $r_G(x,y)=0$ if and only if $x=y$,\\
(3) $r_G(x,y)=r_G(y,x)$,\\
(4) $r_G(x,y)+r_G(y,z)\ge r_G(x,z).$
}

\section{Proof of the Main Results}
In order to determine extremal polygonal chains with respect to Kirchhoff index, we first give some comparison results on Kirchhoff index in these graphs. Recall that we denote the polygons of $P_n$ by $H_1,H_2,\ldots,H_{n}$ such that $H_i$ is adjacent to $H_{i+1}\,(1\leq i\leq n-1),$ and we use $P(w)$ to denote a polygonal chain with $n$ polygons and $w=(w_2,w_3,\ldots,w_{n-1})$ is a $(n-2)$-tuple of $0, 1,\ldots, k_i-4$.

 \begin{figure}[!ht]
\centering
  \includegraphics[width=120mm]{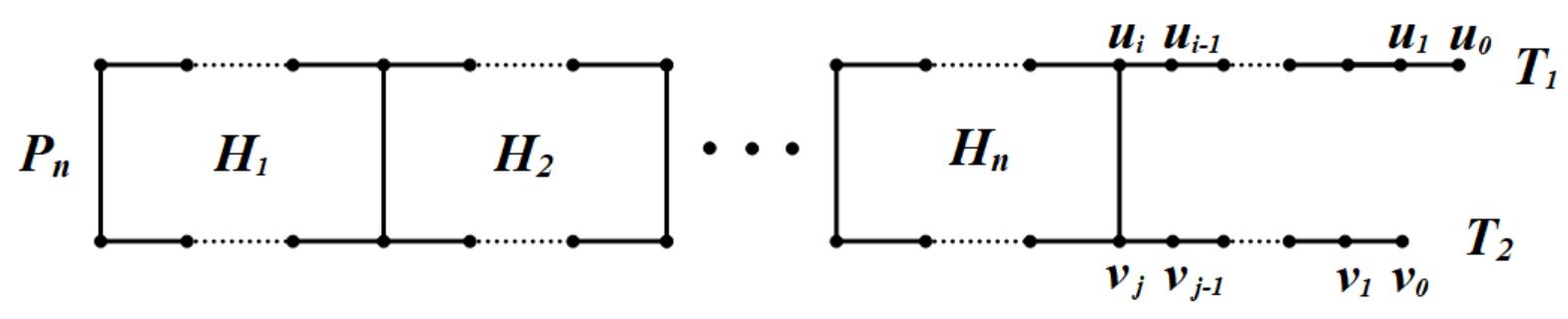}
  \caption{The graph $N$ in the proof of Lemma 8.}
\end{figure}

{\bf Lemma 8.} {\it Let $T_1$ be a $(u_i,u_0)$-path of length $i$ and let $T_2$ be a $(v_j,v_0)$-path of length $j.$ The graph $N$ obtained by attaching the endpoints $u_i$ and $v_j$ of $T_1$ and $T_2$ to two adjacent vertices of degree $2$ in $H_{n}$. If $i>j,$ we have $r_N(u_0)>r_N(v_0).$}

{\bf Proof.} We first denote by $u_i,u_{i-1},\ldots,u_0$ the vertices of $T_1$ and denote by $v_j,u_{j-1},\ldots,v_0$ the vertices of $T_2;$ See Fig. 4. Note that if $k\in V(P_n),$ then $r_{N}(u_0,k)=r_{N}(u_i,k)+i$ and $r_{N}(v_0,k)=r_{N}(v_j,k)+j.$ Observe that $V(N)=V(P_n)\setminus \{u_i,v_j\}\bigcup V(T_1)\bigcup V(T_2),$ we will distinguish two cases.

Case 1. $k\in V(P_n)\setminus \{u_i,v_j\}.$

By Lemma 7, we have
\begin{align}
\begin{split}
r_{N}(k,v_0)\leq& r_{N}(k,u_i)+r_{N}(u_i,v_j)+r_{N}(v_j,v_0)\\
<&r_{N}(k,u_i)+1+j\\
\leq& r_{N}(k,u_i)+i=r_{N}(k,u_0),
\end{split}
\end{align}
where the second inequality follows from the fact that $r_{N}(u_i,v_j)<d_{N}(u_i,v_j)=1.$

Case 2. $k\in V(T_1)\cup V(T_2).$

\begin{align}
\begin{split}
  \sum_{k\in V(T_1)\cup V(T_2)}r_N(v_0,k)&=\sum_{t=1}^jr_N(v_0,v_t)+ \sum_{t=0}^ir_N(v_0,u_t)\\
&=\frac{j(j+1)}{2}+ j(i+1)+(i+1)r_N(u_i,v_j)+\frac{i(i+1)}{2} \\
&<\frac{i(i+1)}{2}+i(j+1)+(j+1)r_N(u_i,v_j)+\frac{j(j+1)}{2}\\
&=\sum_{t=0}^ir_N(u_0,u_t)+\sum_{t=1}^jr_N(u_0,v_t)\\
&=\sum_{k\in V(T_1)\cup V(T_2)}r_N(u_0,k).
\end{split}
\end{align}
Note that
\begin{align}
\begin{split}
&r_{N}(u_0)-r_{N}(v_0)\\
=&\sum_{k\in V(N)}r_{N}(u_0,k)-\sum_{k\in V(N)}r_{N}(v_0,k)\\
=&\sum_{ k\in V(P_n)\setminus \{u_i,v_j\}}r_{N}(u_0,k)-\sum_{ V(P_n)\setminus \{u_i,v_j\}}r_{N}(v_0,k)\\
&+\sum_{ k\in V(T_1)\cup V(T_2)}r_{N}(u_0,k)-\sum_{ k\in V(T_1)\cup V(T_2)}r_{N}(v_0,k).
\end{split}
\end{align}
By (1), (2) and (3), we deduce that $r_{N}(u_0)-r_{N}(v_0)>0.$ This completes the proof of Lemma 8. \hfill $\Box$

 \begin{figure}[!ht]
\centering
  \includegraphics[width=100mm]{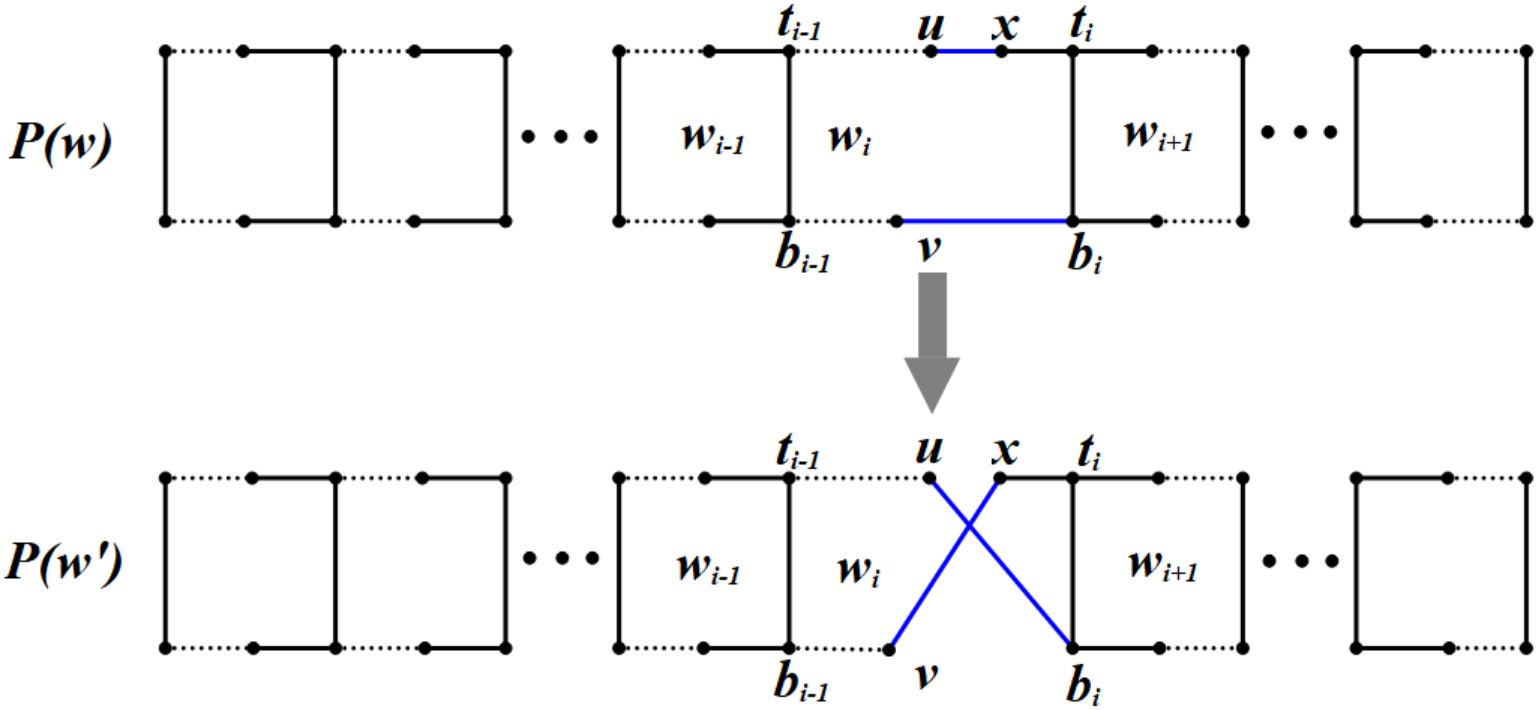}
  \caption{$P(w)$ and $P(w')$ in the proof of Lemma 9.}
\end{figure}

{\bf Lemma 9.} {\it Let $w=(w_2,w_3,\ldots,w_{n-1})$ and $w'=\{w_2,\ldots,w_{i-1},w_i-1,k_{i+1}-4-w_{i+1},\ldots,k_{n-1}-4-w_{n-1}\},$
If there exists $i\in\{2,\ldots,n-1\}$ such that $w_i\ge \frac{k_i-4}{2}+1,$ we have $K\!f(P(w))<K\!f(P(w')).$}

{\bf Proof.} Let $t_ib_i$ be the common edge of $H_i$ and $H_{i+1}$ such that $t_i$ is the top vertex; see Fig 5. Denote by $x$ the neighbour of $t_i$ in $H_i$ and denote by $u$ the other neighbour of $x.$ Let $v$ be the neighbour of $b_i;$ see Fig 5. Note that the edge set $\{ux\,,vb_i\}$ is an edge cut of $P(w).$ Deleting edges $\{ux\,,vb_i\}$ and adding edges $\{ub_i\,,vx\},$ we can get graph $P(w').$ According to the construction of $P(w')$, we deduce that $P(w)$ and $P(w')$ are pairs of $S,T$-isomers. Let $A_1$ be the component of $P(w)-\{ux\,,vb_i\}$ such that $\{u,v\}\subseteq V(A_1)$ and let $B_1$ be the component of $O(w)-\{ux\,,vb_i\}$ such that $\{x,b_i\}\subseteq V(B_1).$ By Lemma 6, We have
$$
K\!f(P(w))-K\!f(P(w'))=\frac{[r_{A_1}(u)-r_{A_1}(v)][r_{B_1}(b_i)-r_{B_1}(x)]}{r_{A_1}(u,v)+r_{B_1}(x,b_i)+2}.
$$
By Lemma 8, we deduce that $r_{A_1}(u)-r_{A_1}(v)>0$ and $r_{B_1}(b_i)-r_{B_1}(x)<0.$ It follows that
$$
K\!f(P(w))-K\!f(P(w'))<0.
$$
This completes the proof of Lemma 9. \hfill $\Box$

Through similar discussions as Lemma 9, we can get the following Lemma 10.

{\bf Lemma 10.} {\it  Let $w=(w_2,w_3,\ldots,w_{n-2})$ and $w'=\{w_2,\ldots,w_{i-1},w_i+1,k_{i+1}-4-w_{i+1},\ldots,k_{n-1}-4-w_{n-1}\},$
If there exists $i\in\{1,2,\ldots,n-1\}$ such that $\frac{k_i-4}{2}-1\ge w_i,$  we have $K\!f(P(w))<K\!f(P(w')).$}

Recall that we denote the polygons of $P_n$ by $H_1,H_2,\ldots,H_{n}$ such that $H_{i}$ is adjacent to $H_{i+1}(1\leq i\leq n-1).$ Moreover, let $t_ib_i$ be the common edge of $H_{i}$ and $H_{i+1}$ such that $t_i$ is the top common vertex.
\begin{figure}[!ht]
\centering
  \includegraphics[width=160mm]{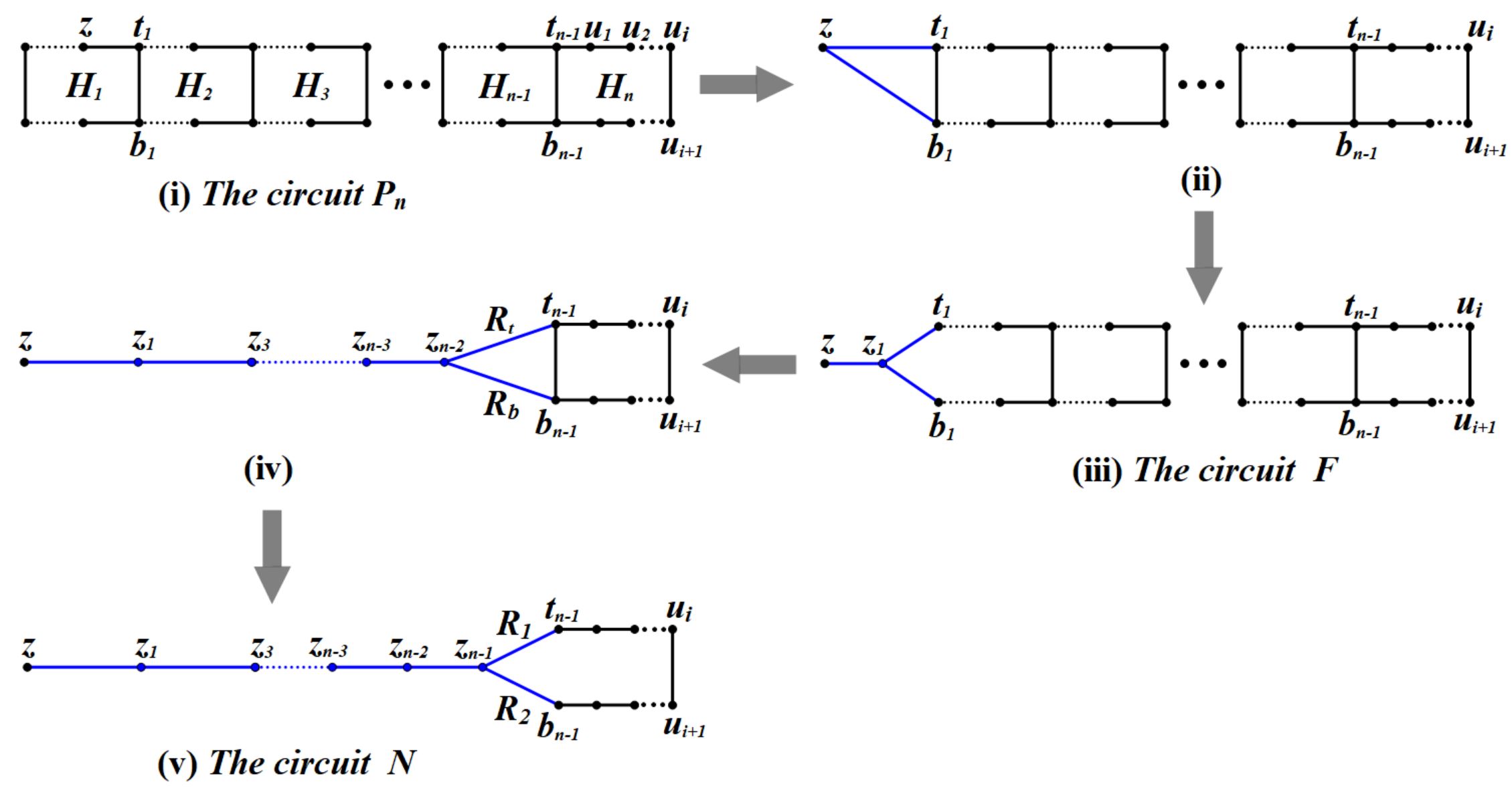}
  \caption{Illustration of circuit simplification to $P_n$ in the proof of Lemma 11.}
\end{figure}

{\bf Lemma 11.} {\it Denote the vertices of $H_n$ by $t_{n-1},u_1,u_2,\ldots,u_{k_n-2},b_{n-1}$ respectively. If $P_n$ is a weighted polygonal chain and the weight on edge $t_{n-1}b_{n-1}$ is $1,$ then for any $z\in V(H_1)\backslash \{t_1,b_1\},$ we have $r_{P_n}(z,u_{i})<r_{P_n}(z,u_{i+1}),$ $1\leq i \leq \lfloor (k_{n}-3)/2\rfloor.$}

{\bf Proof.} Note that $V(H_1)\backslash \{t_1,b_1\}$ has $k_1-2$ different vertices, we first choose the vertex $z$  which is adjacent to $t_1.$ The other remaining vertices can be proved in the same way, so we omit the process. In order to obtain our result, we need the following steps to simplify the circuit $P_n$.

$\bullet$ First, perform the Series transform on $H_1$ of $P_n$ to turn it into a triangle as shown in Fig. 6(ii).

$\bullet$ Next, perform the $\Delta$-$Y$ transform on this new triangle. This results in a new vertex $z_1$ as shown in Fig. 6(iii). We denote this resulting circuit by $F.$

By Lemma 7, we have $$r_{P_n}(z,u_{i})=r_{F}(z,u_{i}),\,\,\,\, r_{P_n}(z,u_{i+1})=r_{F}(z,u_{i+1}).$$

By repeatedly using this two steps, we obtain the simplified circuit of $P_n$ as depicted in Fig. 6(iv). Note that the edges $z_{n-2}t_{n-1}$ and $z_{n-2}b_{n-1}$ in Fig. 6(iv) are the new edges after the transformation. Denote the weighes of $z_{n-2}t_{n-1}$ and $z_{n-2}b_{n-1}$ by $R_t$ and $R_b.$ By making $\Delta$-$Y$ transformation to triangle $z_{n-2}t_{n-1}b_{n-1},$ we could obtain a simplified circuit $N,$ as shown in Fig. 6(v). Suppose that the weights of edges  $z_{n-1}t_{n-1}$ and $z_{n-1}k_{n-1}$ in $N$ are $R_1$ and $R_2$. Recall that the weight of edge $t_{n-1}b_{n-1}$ is $1.$ By Definition 3, we have
\begin{eqnarray*}
  R_1=\frac{R_t}{R_t+R_b+1}, \ \  R_2=\frac{R_b}{R_t+R_b+1}.
\end{eqnarray*}
  Thus $0<R_1<1.$ By using the parallel and series circuit reductions rules, we have
\begin{eqnarray*}
  &r_{P_n}(z,u_{i})=r_N(z,u_{i})=r_N(z,z_{n-1})+\frac{(R_1+i)(R_2+k_n-1-i)}{R_1+R_2+k_n-1}, \\
  &r_{P_n}(z,u_{i+1})=r_N(z,u_{i+1})=r_N(z,z_{n-1})+\frac{(R_1+i+1)(R_2+k_n-2-i)}{R_1+R_2+k_n-1}.
\end{eqnarray*}
Hence, we have
\begin{align*}
  r_{P_n}(z,u_{i})-r_{P_n}(z,u_{i+1})&=\frac{(R_1+i)(R_2+k_n-1-i)}{R_1+R_2+k_n-1}-\frac{(R_1+i+1)(R_2+k_n-2-i)}{R_1+R_2+k_n-1}\\
  &=\frac{R_1-R_2-k_n+2i+2}{R_1+R_2+k_n-1}<0.
\end{align*}
The third inequality follows from the conditions $0<R_1<1$ and $1\leq i \leq \lfloor (k_{n}-3)/2\rfloor.$ Therefore,we have
$$
 r_{P_n}(z,u_{i})<r_{P_n}(z,u_{i+1}).
$$
This completes the proof of Lemma 11. \hfill $\Box$

\begin{figure}[!ht]
\centering
  \includegraphics[width=120mm]{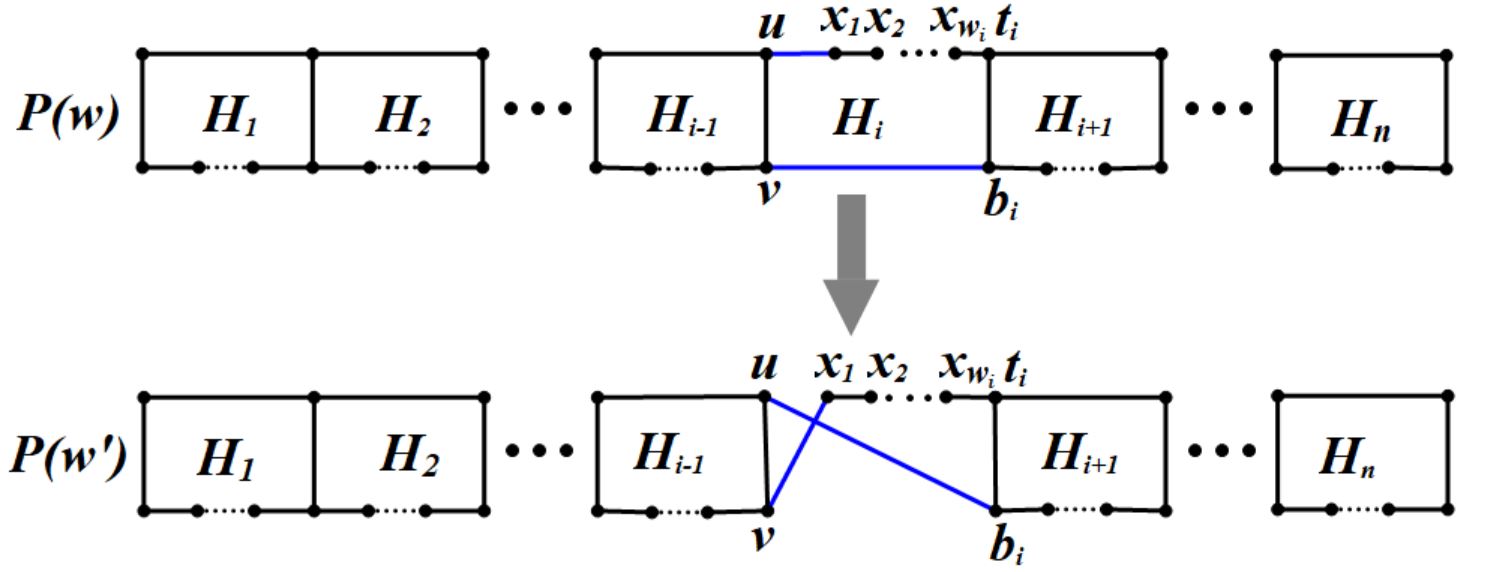}
  \caption{Illustration of circuit simplification to $P(w)$ in the proof of Theorem 1.}
\end{figure}

Now we are ready to proof our main results.

{\bf Proof of Theorem 1.} Let $w=(w_2,w_3,\ldots,w_{n-1}).$ Suppose that $P(w)$ has minimum Kirchhoff index. By Lemma 9 and 10, we have $w_i=0$ or $k_i-4.$ Without loss of generality, we assume that $w_2=0.$ Now we are going to prove $w_i=0$ for all $i\in \{2,\ldots,n-1\}.$ On the country, there exists some $i$ such that $w_{i-1}=0$ and $w_i=k_i-4,$ denote the vertices of the polygon $H_i$ by $u,x_1,\ldots,x_{w_i},t_{i},b_{i},v$ as show in Fig. 7. Let $P(w')$ be the graph obtained from $P(w)$ by deleting edges $\{ux_1,vb_{i}\}$ and adding two new edges $\{ub_{i},vx_1\}.$ Next, we will prove $K\!f(P(w'))<K\!f(P(w)).$

Note that $P(w)$ and $P(w')$ are pairs of $S,T$-isomers. Let $A$ be the component of $P(w)-\{ux_1,vb_{i}\}$ such that $\{u,v\}\subseteq V(A),$ and  Let $B$ be the component of $P(w)-\{ux_1,vb_{i}\}$ such that $\{x_1,b_i\}\subseteq V(B).$ By Lemma 6, we have
\begin{equation}
K\!f(P(w))-K\!f(P(w'))=\frac{[r_A(u)-r_A(v)][r_B(b_{i})-r_B(x_1)]}{r_A(u,v)+r_B(x_1,b_{i})+2}.
\end{equation}

By Lemma 8, we can obtain $r_{B}(b_i)-r_{B}(x_1)<0.$ We then consider $r_A(u)$ and $r_A(v).$ Let $z$ be a vertex in $A,$ we distinguish two cases.

Case 1. $z\in V(H_j),1 \leq j\leq i-2.$ By using series and parallel connection rules, we can simplify $A$ to a weighted polygonal chain consisting of polygons $H_j, H_{j+1},\ldots,H_{i-1}.$ Note that the weight on edge $t_{i}b_{i}$ is $1$. By Lemma 11, we have $r_A(z,u)<r_A(z,v).$

Case 2. $z\in V(H_{i-1}).$ By using series and parallel connection rules, we can simplify $A$ to a weighted polygon $H_{i-1}$ with the weight $r$ on edge $t_{i-2}b_{i-2}$ and the weight $1$ on other edges. Then
\begin{align*}
  &\sum_{z\in V(H_{i-1})}r_A(z,u)
  =\frac{r+k_{i-1}-2}{r+k_{i-1}-1}+\sum_{j=1}^{k_{i-1}-2}\frac{(r+j)(k_{i-1}-1-j)}{r+k_{i-1}-1},\\
  &\sum_{z\in V(H_{i-1})}r_A(z,v)
  =\frac{r+k_{i-1}-2}{r+k_{i-1}-1}+\frac{2(r+k_{i-1}-3)}{r+k_{i-1}-1}
  +\sum_{j=2}^{k_{i-1}-2}\frac{(r+j)(k_{i-1}-1-j)}{r+k_{i-1}-1}.
\end{align*}
Noting that the initially weight of edge $t_{i-1}b_{i-1}$ is $1,$ we have $r<1.$ Since
\begin{equation*}
  \sum_{z\in V(H_{i-1})}r_A(z,u)-\sum_{z\in V(H_{i-1})}r_A(z,v)=\frac{(k_{i-1}-4)(r-1)}{r+k_{i-1}-1}\le0,
\end{equation*}
we have
\begin{equation*}
 \sum_{z\in V(H_{i-1})}r_A(z,u)\le\sum_{z\in V(H_{i-1})}r_A(z,v).
\end{equation*}
By Cases 1-2, we deduce that $r_A(u)<r_A(v).$  Combining with the fact $r_{B}(b_i)-r_{B}(x_1)<0$, we can obtain
$$
K\!f(P(w))>K\!f(P(w')).
$$
A contradiction. This completes the proof of Theorem 1. \hfill $\Box$

{\bf Proof of Theorem 2.} Let $P(w)$ be a polygonal chain with the maximum Kirchhoff index. By Lemma 9 and Lemma 10, we have $|2w_i-(k_i-4)|\le 1.$  So we have $w_i=\lfloor\frac{k_i-4}{2}\rfloor$ or $\lceil\frac{k_i-4}{2}\rceil.$ This completes the proof of Theorem 2. \hfill $\Box$

{\bf Proof of Corollary 3.} Note that if $k_i$ is even, we have $w_i=\lfloor\frac{k_i-4}{2}\rfloor=\lceil\frac{k_i-4}{2}\rceil=\frac{k_i-4}{2}.$
So Corollary 3 can be obtained immediately from Theorem 2.\hfill $\Box$

For two vertices $u$ and $v$, we use the symbol $u\leftrightarrow v$ to mean that $u$ and $v$ are adjacent.

{\bf Proof of Corollary 4.} Let $P(w)$ be the odd polygonal chain with the maximum Kirchhoff index. By Theorem 2, we have $w_i=\frac{k_i-3}{2}$ or $\frac{k_i-5}{2}$ for all $i\in \{2,3,\ldots,n-1\}.$  Suppose that there exists  $w_i=\frac{k_i-5}{2}$ and $w_{i+1}=\frac{k_{i+1}-5}{2}.$ Let $x,y \in V(H_{i+1})$ such that $t_i\leftrightarrow x$ and $b_i\leftrightarrow y.$  Let $P(w')$ be the graph obtained from $P(w)$ by deleting edges $\{t_ix,b_iy\}$ and adding two new edges $\{t_iy,b_ix\}.$ Next, we will prove $K\!f(P(w))<K\!f(P(w')).$

By the construction of $P(w'),$ we have that $P(w)$ and $P(w')$ are pairs of $S,T$-isomers. Let $A$ be the component of $P(w)-\{t_ix,b_iy\}$ such that $\{t_i,b_i\}\subseteq V(A),$ and Let $B$ be the component of $P(w)-\{t_ix,b_iy\}$ such that $\{x,y\}\subseteq V(B).$ By Lemma 6, we have
\begin{equation}
K\!f(P(w))-K\!f(P(w'))=\frac{[r_A(t_i)-r_A(b_i)][r_B(y)-r_B(x)]}{r_A(t_i,b_i)+r_B(y,x)+2}.
\end{equation}

By Lemma 6, we can also obtain $r_{B}(y)-r_{B}(x)>0.$ Then we consider $r_A(t_i)$ and $r_A(b_i).$ Let $z$ be a vertex in $A.$ If $z\in V(A)\setminus V(H_{i}),$ by Lemma 8, we have $r_A(z,t_i)<r_A(z,b_i).$ Otherwise,  by using series and parallel connection rules, we can simplify $A$ to a weighted polygon $H_{i}$ with the weight $r$ on edge $t_{i-1}b_{i-1}$ and the weight $1$ on other edges. Then
\begin{align*}
  &\sum_{z\in V(H_{i})}r_A(z,t_i)
  =\sum_{j=1}^{w_i+1}\frac{j(r+k_i-j-1)}{r+k_i-1}+\sum_{j=0}^{k_i-w_i-3}\frac{(r+w_i+j+1)(k_{i}-w_i-j-2)}{r+k_{i}-1},\\
  &\sum_{z\in V(H_{i})}r_A(z,b_i)
  =\sum_{j=1}^{w_i+2}\frac{j(r+k_i-j-1)}{r+k_i-1}+\sum_{j=1}^{k_i-w_i-3}\frac{(r+w_i+j+1)(k_{i}-w_i-j-2)}{r+k_{i}-1}.
\end{align*}
Noting that the initially weight of edge $t_{i-1}b_{i-1}$ is $1,$ we have $r<1.$ Since
\begin{align*}
  \sum_{z\in V(H_{i})}r_A(z,t_i)-\sum_{z\in V(H_{i})}r_A(z,b_i)=&\frac{(r+w_i+1)(k_i-w_i-1)}{r+k_{i}-1}-\frac{(w_i+2)(r+k_i-w_i-3)}{r+k_{i}-1}\\
=&\frac{(r-1)(k_i-2w_i-4)}{r+k_{i}-1}<0,
\end{align*}
we have
\begin{equation*}
 \sum_{z\in V(H_{i})}r_A(z,t_i)<\sum_{z\in V(H_{i})}r_A(z,b_i).
\end{equation*}
Thus $r_A(t_i)<r_A(b_i).$ Combine with the fact  $r_{B}(y)-r_{B}(x)>0$, we can obtain
$K\!f(P(w))<K\!f(P(w')).$ A contradiction. Therefore, there is no $i$ such that $w_i=\frac{k_i-5}{2}$ and $w_{i+1}=\frac{k_{i+1}-5}{2}.$ By a similar proof, we have that there is no $i$ such that $w_i=\frac{k_i-3}{2}$ and $w_{i+1}=\frac{k_{i+1}-3}{2}.$
This completes the proof of Corollary 4. \hfill $\Box$

{\bf Acknowledgement.}  I am very grateful to my best friend Leilei Zhang for his constant support. This research was supported by the NSFC grant 12271170.

\end{document}